\documentclass[12pt]{article}
\usepackage{amsmath}
\usepackage{amsthm}
\usepackage{amsfonts}
\newcommand{\bR}{{\mathbb R}}
\newcommand{\bT}{{\mathbb T}}
\newcommand{\bZ}{{\mathbb Z}}

\newcommand{\Tt}{{\mathbb T^3}}
\newtheorem{theorem}{Theorem}
\newtheorem{lemma}{Lemma}
\newtheorem{corollary}{Corollary}
\begin{document}

\begin{center}
{\Large \bf
Characterization of the set of ``ergodic directions''
in the Novikov's problem of quasi-electrons orbits in normal metals
}
\\[20pt]
Roberto De Leo\\
Dept. of Mathematics, U. of Cagliari, deleo@unica.it\\
Dept. of Mathematics, U. of Maryland, rdl@math.umd.edu
\end{center}

\vskip .5truecm

\begin{abstract}
\sl Novikov's problem of semiclassical orbits of quasi-electrons
in a normal metal leads to a correspondance between 3-ply periodic
functions in $\bR$ and fractals in $\bR P^2$. These fractals are 
the complement of infinitely many open sets labeled by integer 2-cycles
of $\bT^3$. Here we present a characterization of the fractal points
in terms of the open sets labels.
\end{abstract}
\vskip 1.cm
It was well known by physicists since late fifties that the qualitative behaviour 
of conduction in normal metals under a strong uniform magnetic field is dictated 
by the topological properties of the orbits of electrons quasi-momenta. 
Although, no substantial progress was made from the topological point
of view until eighties, when Novikov [Nov82] noticed that a beautiful
topological structure was hidden inside this problem and its pupils
found its main properties.
\par
In a normal metal, i.e. an ions lattice $\Gamma\subset\bR^3$, quasi-momenta
$(p_1,p_2,p_3)\in\left(\bR^3\right)^*$ are defined modulo a vector of the 
dual lattice $\Gamma^*\subset\left(\bR^3\right)^*$ and their orbits are given by 
the intersection between a level surface of the dispersion law $\varepsilon(p)$
and the bundle of planes perpendiculat to the magnetic field 
$\omega\in\Omega^1(\bT^3)$. 
Therefore the mathematical setting of this Novikov problem is very simple:
given a smooth function $f:\bT^3\to\bR$ we want to study the existence and the 
topological properties of the open 
intersections (if any) between a level surface $M^2_c$ of $f$ and the level sets of a 
constant 1-form $\omega$ as a function of the direction of $\omega$ and of the 
level sets of $f$.
\par
After the fundamental results found by Zorich [Zor84] and Dynnikov [Dyn93,Dyn97] 
the following picture emerged: every function $f$ induces on the space $\bR P^2$ 
of directions of $\omega$ two functions $c_{m,M}:\bR P^2\to\bR$ s.t. the
set $c_m(\omega) \ne c_M(\omega)$ is the disjoint union of open sets $S_i$,
each of them labeled by an element $l_i$ of $H_2(\bT^3,\bZ)$; moreover, the 
complement $F$ of the union of $S_i$ has a fractal structure.
\par
The meaning of all this construction is the following:
if $c_m(\omega)<c<c_M(\omega)$ then $\omega$, that belongs to some $S_{i_0}$, 
induces on $M^2_c$ both closed
and open orbits and such open orbits are strongly asymptotic to the intersection
between a plane perpendicular to $\omega$ and $l_{i_0}$
(in the universal covering), while if $c<c_m(\omega)$ or $c_M(\omega)<c$ only
closed orbits appear. The set $c_m(\omega)=c_M(\omega)$ instead is the union of
the boundaries $\partial S_i$ and of the set $E$ of the so-called ``ergodic
directions'', i.e. the directions of $\omega$ that induce on $M^2_c$ open orbits
that fill components of genus bigger than 2 [DL99].
\par
In this communication we present a simple characterization of the ergodic
direction with maximal degree of irrationality:
\par
\medskip
\begin{theorem}
	The 3-irrational accumulation points of $\{l_i\}$ are exactly the
	3-irrational ergodic directions.
\end{theorem}
\medskip
To prove this theorem we need to prove three properties of the set
$E$ of ergodic directions, that has being proved to be not empty in
[DL99].
\par
\medskip
\begin{lemma} 
For any generic function $f$, in any neighborhood $U_\omega$ of 
any direction $\omega\in\cup \partial S_i\cup E$ there are infinitely 
many stability zones.
\end{lemma}
\begin{proof}
	In both cases the neighborhood will not be contained in a single zone.
	Every rational inside $U_\omega$ (or inside 
	$U_\omega\setminus(S_i\cup\partial S_i)$ if the
	point lies on boundary of zone $i$) is contained in some stability zone
	[Dyn99], so there is at least another class. Removing from $U_\omega$ 
	a finite number of classes that do not contain it cannot exhaust its
	points because the boundaries of these closed sets cannot coincide
	in more than a countable set of points, so there will be always some
	open set left and therefore there will be inifitely many zones.\par
	As the size of a zone must go to zero as $\|l_i\|\to\infty$, infinitely 
	many of them will be completely contained inside it.
\end{proof}
\medskip
The second result we need is a recent Dynnikov's proposition [Dyn99] that claims 
that the direction of every label converges to the corresponding stability zone:
\medskip
\begin{lemma} 
$d(S_i,l_i)\leq O(\|l_i\|^{-1})$ for $\|l_i\|\to\infty$
\end{lemma}
\medskip
Let us point out that for this proposition, essential for our main result, Dynnikov 
offers only convincing arguments and not a rigorous proof, that we plan to provide
soon.
\par
Dynnikov's lemma has an obvious corollary that up to now has never been explicitly 
stated:
\medskip
\begin{corollary}
	Given any function $f$, there are infinitely many directions in 
	$H_2(\Tt,\bZ)$ that do not correspond to any stability zones.
\end{corollary}
\begin{proof}
	Choose any point inside a stability zone and be $d$ its distance from the
	boundary. Close to it, say within $d/2$ from it, there are inifinitely many 
	1-rational directions such that, when represented by elements of 
	$H_2(\Tt,\bZ)$, their length will be so big that the zone with that name 
	will have to lie inside a disc of radius smaller than $d/2$. 
	Of course no such homology class but one can be the label of a stability zone.
\end{proof}
\medskip
Last lemma illustrates the relation between the labels of the stability zones $S_i$ and
the complement of their union $\cup S_i$:
\par
\medskip
\begin{lemma} $\overline{\{l_i\}}=\cup \partial S_i\cup E$ \end{lemma}
\begin{proof}
	It is basically a corollary from Dynnikov's lemma and our lemma 1. Indeed, let 
	$\omega\in\cup \partial S_i\cup E$ and consider the balls $B(\omega,\epsilon)$ and
	$B(\omega,\epsilon^\prime)$: inside both balls, by our lemma, there are 
	inifitely many zones, and, in particular, for every $\epsilon$ we can 
	choose an $\epsilon^\prime\leq\epsilon$ so small that all zones 
	fully contained in it will have a label so big that $d(S_j,l_j)\leq\epsilon/2$.
	Taking now $\epsilon^{\prime\prime}=\min(\epsilon^\prime,\epsilon/2)$
	we have that $S_j\subset B(\omega_i,\epsilon/2)$ and $d(S_j,l_j)\leq\epsilon/2$,
	so that finally $d(\omega_i,l_j)\leq\epsilon$, i.e. $l_j\in B(\omega,\epsilon)$.\par
	Taking $\epsilon=1/n$ we can use the above construction to generate a 
	sequence of elements in $\{l_j\}$ converging to $\omega$.
\end{proof}
\medskip
Putting together this few facts we are able to prove easily the main thorem:
\medskip
\begin{proof}
	Boundaries contain at most 2-irrational directions, as they contain 
	a rank-1 subset of lattice points, and therefore all 3-irrational
	directions asymptotically reached by $l_i$ must be ergodic, and all
	of them are reached because of Lemma 3.
\end{proof}

\noindent{\bf Acknowledgments}\par
The author is in debt with his advisor S.P. Novikov for many fruitful
discussions and suggestions and with I. Dynnikov for many technical 
discussions on the topic. The author also acknowledge financial support 
from Indam for its PhD at the UMD and from the Math Dept. at the
U. of Roma for the project Cofin2000 ``Propriet\`a geometriche delle
variet\`a reali e complesse''

\end{document}